\documentclass[10pt]{article}
\usepackage[utf8]{inputenc}
\usepackage[english]{babel}
\usepackage{amsfonts}
\usepackage{amssymb,enumitem}
\usepackage{amsmath,amsthm}
\usepackage{algorithm}
\usepackage{booktabs}
\usepackage{graphicx}
\usepackage{titlesec}
\usepackage{algorithmic}
\usepackage{adjustbox}
\usepackage{multirow}
\usepackage{soul}
\usepackage{caption}
\usepackage{cite}
\usepackage{url}
\theoremstyle{definition}

\newtheorem{thm}{Theorem}[section]
\newtheorem{lemma}{Lemma}[section]

\newtheorem{alg}{Algorithm}[section]
\newcommand{\R}{\mathbb R}

\newcommand{\V}{\mathcal V}
\newcommand{\Rn}{\mathbb R^{n}}

\newcommand{\N}[1]{\mathcal N_{#1}}
\newcommand{\fsub}[2]{f_{\mathcal N_{#1}}(x_{#2})}
\newcommand{\gfsub}[2]{\nabla f_{\mathcal N_{#1}}(x_{#2})}
\newcommand{\I}[1]{\mathcal I_{#1}}

\newcommand{\gfsubi}[2]{\nabla f_{\mathcal I_{#1}}(x_{#2})}

\newcommand{\ke}{k_\varepsilon}

\begin {document}
\title{Subsampled Nonmonotone Spectral Gradient Methods}

\author{
Stefania Bellavia
 \thanks{ Department of Industrial Engineering, University of Florence, Viale Morgagni, 40/44, 50134 Florence, Italy, e-mail: {\tt
stefania.bellavia@unifi.it}. Research supported by  Gruppo Nazionale per il Calcolo Scientifico,
(GNCS-INdAM) of Italy. }
 \and Nata\v sa Krklec Jerinki\'c \thanks{Department of Mathematics and
Informatics, Faculty of Sciences,  University of Novi Sad, Trg
Dositeja Obradovi\'ca 4, 21000 Novi Sad, Serbia, e-mail: {\tt
natasa.krklec@dmi.uns.ac.rs}. Research supported by Serbian Ministry of Education, Science and Technological Development, grant no. 174030.}
\and Greta Malaspina\thanks{Department of Mathematics and
Informatics, Faculty of Sciences,  University of Novi Sad, Trg
Dositeja Obradovi\'ca 4, 21000 Novi Sad, Serbia, e-mail: {\tt
greta.malaspina@dmi.uns.ac.rs}}}
\date{}
\maketitle

\begin{abstract}

\noindent This paper deals with subsampled spectral gradient methods for minimizing finite sums.
Subsample function and gradient approximations are employed  in order to reduce the overall computational cost of  the classical spectral gradient methods. The global convergence is enforced by a nonmonotone line search procedure.  Global convergence is proved provided that functions and gradients are approximated with increasing accuracy.
R-linear convergence  and worst-case iteration complexity is investigated
in case of strongly convex objective function.
Numerical results  on well known binary classification problems are given  to show the effectiveness of
this framework and analyze the effect of different spectral coefficient approximations arising from the variable sample nature of this procedure.


\end{abstract}

\section{Introduction}
The aim of this work is to present a class of first-order iterative methods for optimization problems where the objective function is given as the
mean of a large number of functions, i.e.,
\begin{equation}\label{sum_problem} \min_{x\in \Rn}f_{\N{}}(x), \phantom{spazio}f_{\N{}}(x)=\frac{1}{N}\sum_{j=1}^Nf_j(x).\end{equation}
The functions $f_j:\Rn\rightarrow R$ are assumed to be continuously differentiable for $j=1,\ldots,N$  and $\N{}$ denotes the set of indices $\{1,\ldots,N\}$.
The motivation for studying  problems of this form comes from machine learning.  Indeed,
 the training phase of a neural network requires to solve  problems of the form \eqref{sum_problem} where
the number $N$ of functions is generally large enough to prevent the employment
of classical optimization methods, thus leading to the necessity of developing different strategies.
One possible approach is to employ a so-called mini-batch or subsampling strategy \cite{bgm, bkkj, bbn_2017, ebnkjm,
blatt:incrementalgradient, bbn, nocedal:OMLSML, bottou:stochasticgradient,noc2,
Friedlander, johnson:acc, kkj:lemma, mit,gunter:stochasticquasinewton,tan:SG_BB,yang:minibatch_BB}.
That is, an iterative optimization method is applied but instead of considering the whole sum \eqref{sum_problem},
at the beginning of every iteration functions and/or gradient and/or Hessian are approximated using only a subset of the functions $\{f_1,\dots, f_N\}$.

In this paper we focus on variable sample variants of spectral gradient methods with nonmonotone line search. Spectral gradient methods, originally proposed in
\cite{barzilai:spectralgradient} for the solution of unconstrained optimization problems have been widely used and developed
(see \cite{birgin, dai,fletcher, diserafino, raydan:spectralgradient} and references therein).
In these approaches the steplength is adaptively  chosen
and they showed very good practical performance.
Here, we combine the spectral gradient method with a nonmonotone line search and a variable sample strategy with a twofold aim:  retain robustness
and adaptive steplength choice of spectral gradient methods and reduce the overall computational cost of solving  \eqref{sum_problem} by approximating functions and gradients with increasing accuracy.

Given the  sample size, the subsample is randomly generated from $\{1,\dots,N\}$ and we consider two kinds of subsampling.
As we will see, different options arise for the subsampling approximation of  the relevant quantities characterizing  this class of 
methods. These options  give rise to four variants of the subsampling spectral gradient method with line search and 
we analyze the global
convergence of the unifying framework they belong to.  We remark that, thanks to the globalization strategy, the methods show
global convergence and do not
require to choose the steplength by trial, which can
be time consuming in practice.  Focusing on the strongly convex case we prove R-linear convergence of the generated sequence to the minimizer of \eqref{sum_problem}. We also provide iteration complexity of the methods and show that the complexity bound of the corresponding exact method is retained, despite inaccuracy in functions and gradients.

Finally, we present some numerical results on binary classification problems, showing the effectiveness of our approach.

Spectral gradient  methods for problem \eqref{sum_problem} have been investigated in \cite{tan:SG_BB,yang:minibatch_BB}.
In \cite{tan:SG_BB} they are used in combination with the stochastic gradient method, while in \cite{yang:minibatch_BB}
a mini-batch strategy is employed. In both papers
convergence is proved assuming to employ   the full gradient every $m$ iterations (at each outer iteration). 

We differ from this latter approach as we embed the  variable sample spectral method in a  nonmonotone
line search strategy using approximate functions and gradients. In the convergence  analysis  we have  also to take into account inaccuracy  in function values and not only in gradient,
while in \cite{tan:SG_BB,yang:minibatch_BB} only inaccuracy in gradients needs to be handled as the function values are unused. On the other hand,
the employment of the line search procedure allow us to obtain global convergence to a stationary point
irrespectively of convexity of the objective function, differently from \cite{tan:SG_BB,yang:minibatch_BB}.  Nonconvex problems arise in neural networks training, so this relaxation can be of great interest.

  Taking different samples and/or sample sizes in different iterations affects the spectral coefficients, so, besides the convergence analysis, the main goal of this paper is to investigate the performance of different choices of spectral coefficient calculations as they affect the cost per iteration as well.    Notice that
in \cite{tan:SG_BB} the spectral coefficient is
 computed at the outer iterations using the full gradient and it is taken constant along the inner iterations. In \cite{yang:minibatch_BB} it is updated also in the inner iterations using a subsample of the dataset whose dimension depends on the condition number of the problem.  In \cite{SPG-NKNKJ} spectral gradient methods for problems with objective function given in form of mathematical
expectation have been analyzed, however the effect of the choice of the spectral coefficient has not been investigated.

This paper is structured as follows.
In Section 1 we recall the classical nonmonotone spectral gradient method, introduce the subsampling strategy and embed it in the framework
of the classical method. In Section 2 we study the convergence behavior of the obtained subsampling procedure providing convergence and complexity results.  In Section 3 we report the results of numerical tests we carried out. We focus on the comparison
between the subsampled method and the classical counterpart and on the influence of some critical parameters over the performance of the methods.

\section{The Method}
In this section we describe the subsampled spectral gradient framework we are involved with.  Spectral gradient methods are first-order iterative procedures employing
the gradient vector as a
search direction.  A crucial role is played by the steplength that is chosen
 in such a way to inject some second-order information into the methods.

At a generic iteration, given the current point $x_k$, the new iterate is computed as
$$x_{k+1}=x_k+\alpha_kd_k$$
where
\begin{equation}\label{eq:SGdir}
d_k=-\sigma_{k-1}^{-1}\nabla f_{\N{}}(x_k)
\end{equation}
and 
 $\alpha_k$ is chosen through a line search strategy. 
The scalar $\sigma_k$ is called \emph{spectral coefficient}.\\
$\ $\\

Here we adopt the classical Barzilai-Borwein choice given in \cite{barzilai:spectralgradient}, i.e.\\
\begin{equation*}\label{eq:SpecCoeff}
\sigma_k=
\begin{cases}
\displaystyle\frac{s_{k}^ty_{k}}{s_{k}^ts_{k}} & \text{if }\ \displaystyle\frac{s_{k}^ty_{k}}{s_{k}^ts_{k}}\in[\sigma_m, \sigma_M]\\
1 & \text{otherwise}
\end{cases}
\end{equation*}
with
\begin{equation*}
\begin{aligned}
&s_{k}=x_{k+1}-x_{k}\phantom{spazio}\\
&y_{k}=\nabla f_{\N{}}(x_{k+1})-\nabla f_{\N{}}(x_{k})
\end{aligned}
\end{equation*}
and $0<\sigma_m<1<\sigma_M<+\infty$  given safeguards.  Notice that this guaranties that the search direction is a descent direction  with respect to the current objective function. Therefore, provided that the function is bounded from below on the line segment $[x_k, x_k+d_k]$, there exists a steplength interval that satisfies Wolfe conditions.  In practice, Wolfe conditions are often replaced by the backtracking technique combined with the first Wolfe condition, i.e. the Armijo condition. Nonmonotone line search in general allows  even larger step sizes, so the line search remains well defined.  

Specifically, we here assume that the steplength $\alpha_k$ satisfies the Li-Fukushima \cite{lifukushima:ls} nonmonotone descent condition
and the Wolfe condition, i.e. it satisfies the following two inequalities:
\begin{equation}\label{eq:Li}
\begin{aligned}
 &f_{\N{}}(x_k+\alpha_kd_k)\leq f_{\N{}}(x_k)+c_1\alpha_k\nabla f_{\N{}}(x_k)^td_k+\zeta_k\\
 &\nabla f_{\N{}}(x_{k}+\alpha_kd_k)^td_k\geq c_2\nabla f_{\N{}}(x_k)^td_k
 \end{aligned}
\end{equation}
where $0< c_1\leq c_2< 1$, $\zeta_k\geq0$ such that $\sum_{k\geq0}\zeta_k<+\infty.$ 

The following algorithm summarizes a generic iteration of the spectral gradient method we are considering.

\begin{alg}\label{alg:SG}($k$-th iteration of spectral gradient method)$\ $\\
\vspace{0.3cm}
\textbf{Input:} $x_k\in\Rn$, $f:\Rn\rightarrow\R$, $\zeta_k>0$, $0<\sigma_m<1<\sigma_M$\\
\begin{tabular}{p{0.2cm}|p{10cm}}
&set $g_k=\nabla f_{\N{}}(x_k)$\\
&set $s_{k-1}=x_k-x_{k-1}$\\
&set $y_{k-1}=g_k-g_{k-1}$\\
&set $\sigma_{k-1}=\displaystyle\frac{s_{k-1}^ty_{k-1}}{s_{k-1}^ts_{k-1}}$\\
&\textbf{if} $\ \sigma_{k-1}\notin[\sigma_m, \sigma_M]$\\
    & \begin{minipage}{\textwidth}
    \begin{tabular}{c|p{10cm}}
    &set $\sigma_{k-1}=1$
    \end{tabular}
    \end{minipage}\\
  &\textbf{end if}\\
&set $d_k=-\sigma_{k-1}^{-1}g_k$\\
&compute $\alpha_k$ such that \eqref{eq:Li} holds\\
&set $x_{k+1}=x_k+\alpha_kd_k$
\end{tabular}
\vspace{0.3cm}
$\ $\\
\end{alg}

In order to reduce the overall computational cost of the procedure, at each iteration we choose a subsample
$\N{k}\subseteq\{1,\dots,N\}$ of size $N_k$ and we employ a variable sample strategy. That is, we approximate the objective function and its
gradient as follows:
$$ f_{\N{k}}(x):=\frac{1}{N_k}\sum_{j\in\N{k}}f_j(x) $$
$$ \nabla f_{\N{k}}(x):=\frac{1}{N_k}\sum_{j\in\N{k}}\nabla f_j(x). $$

Given an increasing sequence $\{N_k\}$ of sample sizes, we consider two different kinds of subsampling:
\begin{itemize}
 \item \emph{nested subsamples} $\N{k-1}\subseteq\N{k}$;\\
 we take $\N{k}$ as the union of $\N{k-1}$ and a set of $(N_k-N_{k-1})$ randomly-chosen indices in $\mathcal N\setminus\N{k-1}$;
 \item \emph{non-nested subsamples} such that $\N{k-1}\cap\N{k}\neq\emptyset$;\\
 we take one index $j_1$ randomly chosen in $\N{k-1}$ to assure a non-empty intersection, then we take $\N{k}$ as the union of $j_1$
 and a set of $(N_k-1)$ randomly chosen indices in $\N{}\setminus\{j_1\}.$ 
 Notice that we enforce non empty intersection between consecutive subsamples, but since we randomly choose the
indices in $\N{}\setminus\{j_1\}$, we don't have any control over the actual size of the intersection. 

\end{itemize}

In the definition \eqref{eq:SGdir} of the direction $d_k$ we replace the gradient $\nabla f_{\N{}}(x_k)$
with $\gfsub{k}{k}$,
while for the gradient displacement vector $y_{k-1}$
we have different possible choices that we are now going to present.\\
First, let us assume that we are in the nested case (i.e. $\N{k-1}\subseteq\N{k}$).
Since the subsample that we are considering at the current iteration is $\N{k}$, the first option is to replace in the computation of $y_k$ both $\nabla f_{\N{}}(x_k)$ and $\nabla f_{\N{}}(x_{k-1})$
with the approximation given by the current subsample, getting
\begin{equation}\label{eq:Y1}
 y_{k-1}^{(1)}:=\gfsub{k}{k}-\gfsub{k}{k-1}.
\end{equation}
On the other hand we already approximated $\nabla f_{\N{}}(x_{k-1})$ at the previous iteration as $\gfsub{k-1}{k-1}$, therefore the second option is to use 
\begin{equation}\label{eq:Y2}
 y_{k-1}^{(2)}:=\gfsub{k}{k}-\gfsub{k-1}{k-1}.
\end{equation}
Assuming we are in the non-nested case, the definition \eqref{eq:Y1} is still possible, and denoting with $\I{k}$ the intersection of the current and the previous subsample
we can also take
\begin{equation}\label{eq:Y3}
 y_{k-1}^{(3)}:=\gfsubi{k}{k}-\gfsubi{k}{k-1}.
\end{equation}
Each of these choices allows us to exploit a different amount of information in the approximation of $y_{k-1}$ and therefore of the spectral
coefficient $\sigma_{k-1}$, but it also requires a different amount of computation in terms of number of component gradient
evaluations.  We already noticed that at every iteration we need to evaluate
$\gfsub{k}{k}$, then at iteration $k-1$ we have evaluated $\nabla f_{j}(x_{k-1})$ for every $j\in\N{k-1}$.
In the nested case one can store only the previous (average) gradient and the previous sample size and 
the computation of $y_{k-1}^{(2)}$ does not require  extra computation, while  the computation of $y_{k-1}^{(1)}$ requires $(N_k-N_{k-1})$ new evaluations of component gradients.
In the non-nested case,  the storage of all the component  gradients from the previous iteration is needed to compute 
$y_{k-1}^{(1)}$ and   $y_{k-1}^{(3)}$.  The evaluation of  $y_{k-1}^{(1)}$ also requires $(N_k-|\I{k}|)$ evaluations of  the component gradients that have not been already evaluated at the previous iteration.
  If  the storage of all the component gradients
  is not feasible, one can   use a fixed mini-batch for the intersection and use it only for obtaining the spectral coefficient.

The employment of the subsampling scheme to the line search conditions \eqref{eq:Li} is immediate and leads to:
\begin{equation}\label{eq:LiSub}
 f_{\N{k}}(x_k+\alpha_kd_k)\leq f_{\N{k}}(x_k)+c_1\alpha_k\gfsub{k}{k}^td_k+\zeta_k
 \end{equation}
\begin{equation}\label{eq:WolfeSub}
 \nabla f_{\N{k}}(x_{k}+\alpha_kd_k)^td_k\geq c_2\gfsub{k}{k}^td_k
\end{equation}
with all the requests over $c_1, c_2$ and $\{\zeta_k\}$ unchanged.\\
In Table \ref{tabellametodi} we summarize all the possible combinations for the definitions of the vector $y_k$ together
with the two kinds of subsampling, then in Algorithm \ref{alg:MBSM} we present the general structure of the $k$-th iteration
of these methods. The input variable \emph{NM} denotes the name of the method, according to Table \ref{tabellametodi}.

\

\begin{alg}(Framework of variable sample spectral methods, $k$-th iteration)\label{alg:MBSM}$\ $\\
\vspace{0.3cm}
\textbf{Input:} $x_k\in\Rn$, $\{f_j\}_{j\in\N{}}$, $N_0$ initial size, \emph{NM}, $\zeta_k>0$, $0<\sigma_m<1<\sigma_M$\\
\begin{tabular}{p{0.2cm}|p{10cm}}
\footnotesize{1}&{ set $N_k\geq N_{k-1}$}\\
\footnotesize{2}&generate $\N{k}$ according to \emph{NM} and the second column of Table \ref{tabellametodi} \\
\footnotesize{3}&set $g_k=\gfsub{k}{k}$\\
\footnotesize{4}&set $s_{k-1}=x_{k}-x_{k-1}$\\
\footnotesize{5}& compute $y_{k-1}$ according to \emph{NM} and the third column of Table \ref{tabellametodi}\\
\footnotesize{6}&set $\sigma_{k-1}=\displaystyle\frac{s_{k-1}^ty_{k-1}}{s_{k-1}^ts_{k-1}}$\\
\footnotesize{7}&\textbf{if} $\ \sigma_{k-1}\notin[\sigma_m, \sigma_M]$\\
    & \begin{minipage}{\textwidth}
    \begin{tabular}{c|p{10cm}}
    &set $\sigma_{k-1}=1$
    \end{tabular}
    \end{minipage}\\
  &\textbf{end if}\\
\footnotesize{8}&compute $d_k=-\sigma_{k-1}^{-1}g_k$ \\
\footnotesize{9}& compute $\alpha_k$ such that {\eqref{eq:LiSub} and \eqref{eq:WolfeSub} hold}\\
\footnotesize{10}&set $x_{k+1}=\displaystyle x_k+\alpha_kd_k$\\
\end{tabular}
\vspace{0.3cm}
$\ $\\
\end{alg}

\begin{table}[h]
\begin{center}
\caption{Subsampled spectral methods}
{\begin{tabular}{c|c|c}
\hline
 $\phantom{sp}$\textbf{Name}$\phantom{sp}$&  \textbf{Subsample}&$\phantom{spaz}y_k\phantom{spaz}$\\
 \hline
 SG\_N\_1 & Nested   & $y^{(1)}_{k}$ given in \eqref{eq:Y1}\\
 SG\_N\_2 & Nested   & $y^{(2)}_{k}$ given in \eqref{eq:Y2}\\
 SG\_I\_1 & Non-Nested   & $y^{(1)}_{k}$ given in \eqref{eq:Y1}\\
 SG\_I\_3 & Non-Nested   & $y^{(3)}_{k}$ given in \eqref{eq:Y3}\\
  \hline
 \end{tabular}}
\label{tabellametodi}
\end{center}
\end{table}


\section{Global Convergence}
We assume that the objective function $f_{\N{}}$ is bounded from below and
continuously differentiable in $\Rn$, and that each gradient $\nabla f_j$ is Lipschitz-continuous.
We define the errors of approximation $\nu_k$ and $\eta_k$  as follows:
\begin{equation}\label{approxnu}
 \nu_k:=\max\{|\fsub{k}{k}-f_{\N{}}(x_k)|, |\fsub{k}{k+1}-f_{\N{}}(x_{k+1})|\}
\end{equation}
\begin{equation}\label{approxeta}
 \eta_k:=\max\Big\{\big|\|\gfsub{}{k}\|^2-\|\nabla f_{\N{k}}(x_k)\|^2\big|, \big|\|\gfsub{}{k+1}\|^2-\|\nabla f_{\N{k}}(x_{k+1})\|^2\big|\Big\}.
\end{equation}

We prove here that any limit point of the sequence generated by the method is a stationary point of \eqref{sum_problem}
and that when the objective function is strongly convex, $R$-linear convergence to the solution holds.
\begin{thm}\label{thm:statpt}
Let $\{x_k\}$ be the sequence of the iterates computed by Algorithm \ref{alg:MBSM} with
$\sum_{k\geq0}\zeta_k=\bar\zeta<+\infty$. Let us assume that  $\sum_{k\geq0}\nu_k=\bar\nu<+\infty$
and that $\eta_k$ goes to 0 as $k$ goes to $+\infty.$
If $f_{\N{}}$ is bounded from below over $\Rn$ and continuously-differentiable, and the gradients $\nabla f_{j}$
are Lipschitz-continuous with Lipschitz constant $L$,
then:
\begin{equation*}
 \lim_{k\rightarrow+\infty}\|\nabla f_{\N{}}(x_k)\|=0.
\end{equation*}

\begin{proof}
 Subtracting $\nabla f_{\N{k}}(x_k)^td_k$ from both sides of \eqref{eq:WolfeSub}
 and using the Lipschitz-continuity of the gradient we get
  \begin{equation*}\begin{aligned}
  &(c_2-1)\nabla f_{\N{k}}(x_k)^td_k\leq(\nabla f_{\N{k}}(x_k+\alpha_kd_k)-\nabla f_{\N{k}}(x_k))^td_k\leq\\
  &\leq\|\nabla f_{\N{k}}(x_k+\alpha_kd_k)-\nabla f_{\N{k}}(x_k))\|\|d_k\|\leq L\alpha_k\|d_k\|^2.
 \end{aligned}\end{equation*}
 By the definition of $d_k$
 \begin{equation}\label{scalarp}\nabla f_{\N{k}}(x_k)^td_k=-\sigma_{k-1}^{-1}\|\nabla f_{\N{k}}(x_k)\|^2\end{equation}
 hence, from the previous inequality,
 \begin{equation*}
  \alpha_k\geq-\frac{\sigma_{k-1}(c_2-1)}{L}.
 \end{equation*}
 Let us define $C:=-\frac{c_1(c_2-1)}{L}$,
 from \eqref{eq:LiSub}, \eqref{scalarp} and the last inequality we have
 \begin{eqnarray}
 \fsub{k}{k+1}&=&f_{\N{k}}(x_k+\alpha_kd_k)\nonumber \\
                     &\leq& f_{\N{k}}(x_k)+
 \frac{c_1(c_2-1)}{L}\|\nabla f_{\N{k}}(x_k)\|^2+\zeta_k\\\label{eq:optgapthk}
 &\leq& f_{\N{k}}(x_k)-C\|\nabla f_{\N{k}}(x_k)\|^2+\zeta_k.\nonumber
 \end{eqnarray}

 From this inequality and \eqref{approxnu} we have
 \begin{equation}\label{eq:optgapth}
  \begin{aligned}
   &f_{\N{}}(x_{k+1})\leq\fsub{k}{k+1}+\nu_k\leq\\
   &\leq\fsub{k}{k}-C\|\gfsub{k}{k}\|^2+\zeta_k+\nu_k\leq\\
   &\leq f_{\N{}}(x_k)-C\|\gfsub{k}{k}\|^2+\zeta_k+2\nu_k\leq\\
   &\leq f_{\N{}}(x_0)-C\sum_{j=0}^k\|\gfsub{j}{j}\|^2+\sum_{j=0}^k(\zeta_j+2\nu_j),
  \end{aligned}
 \end{equation}
thus
  \begin{equation*}
  \sum_{j=0}^k\|\gfsub{j}{j}\|^2\leq\frac{f_{\N{}}(x_0)-f_{\N{}}(x_{k+1})}{C}+\frac{1}{C}\sum_{j=0}^k(\zeta_j+2\nu_j)
 \end{equation*}
 and taking the limit for $k\rightarrow+\infty$ we get
 \begin{equation*}
  \sum_{k\geq0}\|\gfsub{k}{k}\|^2\leq\frac{f_{\N{}}(x_0)-\lim_{k\rightarrow \infty}f_{\N{}}(x_{k+1})}{C}+\frac{1}{C}(\bar\zeta+2\bar\nu).
 \end{equation*}
 Since $f_{\N{}}$ is bounded from below we get
  \begin{equation}\label{serieconv}
  \sum_{k\geq0}\|\gfsub{k}{k}\|^2<+\infty.
 \end{equation}
 By \eqref{approxeta} we thus have
 $$\lim_{k\rightarrow+\infty}\|\nabla f_{\N{}}(x_k)\|^2\leq \lim_{k\rightarrow+\infty}\left(\|\gfsub{k}{k}\|^2+\eta_k\right)=0$$
 and hence the thesis follows.\\
\end{proof}

\end{thm}
In the next theorem we will show that our procedure needs at most $O(\varepsilon^{-2})$ iterations to provide
$\|\nabla f_{\N{k}}(x_k)\|\leq\varepsilon$. We note that the worst case iteration complexity is of the same order of that of
nonmonotone spectral
gradient methods shown in \cite{Grapiglia:complexity}.

\begin{thm}
Let $\{x_k\}$ be the sequence of the iterates computed by Algorithm \ref{alg:MBSM} with
$\sum_{k\geq0}\zeta_k=\bar\zeta<+\infty.$
If we assume that $\sum_{k\geq0}\nu_k=\bar\nu<+\infty$ and that the gradients $\nabla f_j$ are Lipschitz-continuous with constant $L$,
then for any given $\varepsilon>0$ the generated sequence satisfies
$$
\|\nabla f_{\N{k}}(x_k)\|\leq\varepsilon
$$
in at most
$$
\ke=\left\lceil(f_{\N{}}(x_0)-f_{\N{}}(x_*)+\bar\zeta+2\bar\nu)C^{-1}\varepsilon^{-2}\right\rceil
$$
iterations, where the constant $C$ is given by $C:=-\frac{c_1(c_2-1)}{L}$.

 \begin{proof}
  Let us denote with $\ke$ the first iteration such that $\|\nabla f_{\N{\ke}}(x_{\ke})\|<\varepsilon$.\\
  By \eqref{eq:optgapthk} we get, for every $k<\ke$
  \begin{equation}
   f_{\N{k}}(x_k)-f_{\N{k}}(x_{k+1})+\zeta_k\geq C\|\nabla f_{\N{k}}(x_k)\|^2\geq C\varepsilon^2
  \end{equation}
and by \eqref{approxnu}
 \begin{equation}
   f_{\N{}}(x_k)-f_{\N{}}(x_{k+1})+\zeta_k+2\nu_k\geq C\varepsilon^2.
  \end{equation}
  Taking the sum for $k$ from 0 to $\ke-1$ we get
  \begin{equation}\begin{aligned}
\ke C\varepsilon^2&\leq\sum_{k=0}^{\ke-1}\left(f_{\N{}}(x_k)-f_{\N{}}(x_{k+1})+\zeta_k+2\nu_k\right)\leq\\
  &\leq f_{\N{}}(x_0)-f_{\N{}}(x_{\ke})+\bar\zeta+2\bar\nu
                  \end{aligned}
  \end{equation}
and hence the thesis.
  \end{proof}

\end{thm}

In the sequel we will make use of the next two lemmas.
\begin{lemma}\label{lemma:conv}\cite{nocedal:OMLSML}
 If $f_{\N{}}$ is a continuously differentiable function from $\Rn$ to $\R$, strongly convex with constant $c$, then $f_{\N{}}$ has an unique minimizer $x_*$ and,
 for every $x\in\Rn$, the following inequality holds:
$$2c\left(f_{\N{}}(x)-f_{\N{}}(x_*)\right)\leq\|\nabla f_{\N{}}(x)\|^2.$$
\end{lemma}
\begin{lemma}\cite{kkj:lemma}\label{lemmaRlin}
 If a sequence $\{a_k\}$ converges to 0 $R$-linearly then, for every $\rho\in(0,1)$
 the sequence $\{A_k\}$ given by
 $$A_k:=\sum_{j=0}^k\rho^{j}a_{k-j}$$
 converges to 0 $R$-linearly.
\end{lemma}
We now focus on the strongly convex case. {Denoting with $x_*$ the unique minimizer of $f_{\N{}}$, we}
prove that the optimality gap $\fsub{}{k}-\fsub{}{*}$ tends to zero $R$- linearly,  and consequently the method
drives the optimality gap  below $\varepsilon$, even if approximated gradients and functions are used.

\begin{thm}\label{thm:Rlin}
Let $\{x_k\}$ be the sequence of the iterates computed by Algorithm \ref{alg:MBSM}.
Assume that $f_{\N{}}$ is a strongly convex function from $\Rn$ to $\R$ and that the gradients $\nabla f_j$
are Lipschitz-continuous with constant $L$.
Then:
\begin{enumerate}
 \item[(i)] There exists a constant $\rho\in(0,1)$ such that for every index $k$ the optimality gap satisfies
 \begin{equation*}\begin{aligned}
\fsub{}{k+1}-\fsub{}{*}&\leq\rho^{k+1}\left(\fsub{}{0}-\fsub{}{*}\right)+
\sum_{j=0}^k\rho^{j}\zeta_{k-j}+\\&+\sum_{j=0}^k\rho^{j}(2\nu_{k-j}+C\eta_{k-j})
\end{aligned}\end{equation*}
where $C=-\frac{c_1(c_2-1)}{L}$.
\item[(ii)] If the three sequences $\{\zeta_k\},\ \{\nu_k\},\ \{\eta_k\}$ tend to 0 $\ R$-linearly as $k$ goes
to $+\infty$, then  $\fsub{}{k+1}-\fsub{}{*}$ converges to 0 $R$-linearly.
\end{enumerate}
\begin{proof}$\ $\\
By \eqref{eq:optgapth} and \eqref{approxeta} we have
\begin{eqnarray}
 \fsub{}{k+1}-\fsub{}{*}
   &\leq&\fsub{}{k}-\fsub{}{*}-C\|\gfsub{k}{k}\|^2+
               \zeta_k+2\nu_k \label{optimalgap}\\
   &\leq&\fsub{}{k}-\fsub{}{*}-C\|\gfsub{}{k}\|^2+\zeta_k+2\nu_k+C\eta_k. \nonumber
   \end{eqnarray}
 By Lemma \ref{lemma:conv} we have that
 \begin{equation*}
  \|\gfsub{}{k}\|^2\geq\gamma(\fsub{}{k}-\fsub{}{*})
 \end{equation*}
for every $0<\gamma<\min\{2c, L\}$. Note that $(1-C\gamma)\in (0,1)$. Thus, from \eqref{optimalgap}, we get
 \begin{equation*}
\fsub{}{k+1}-\fsub{}{*}\leq(1-C\gamma)(\fsub{}{k}-\fsub{}{*})+\zeta_k+2\nu_k+C\eta_k.
 \end{equation*}
Iteratively applying this inequality and denoting with $\rho$ the quantity $(1-C\gamma)$, we get
\begin{equation*}\begin{aligned}
\fsub{}{k+1}-\fsub{}{*}&\leq\rho^{k+1}\left(\fsub{}{0}-\fsub{}{*}\right)+
\sum_{j=0}^k\rho^{j}\zeta_{k-j}+\\&+\sum_{j=0}^k\rho^{j}(2\nu_{k-j}+C\eta_{k-j})
\end{aligned}\end{equation*}
and thus (i) holds.\\

Let us define $\omega_{j}:=2\nu_{j}+C\eta_{j}$. By the $R$-linear convergence of $\{\nu_k\}$ and $\{\eta_k\}$ we have that
$\omega_k$ converges to 0 $R$-linearly and thus
%
%
by inequality (i) and Lemma \ref{lemmaRlin} we have (ii).\\
\end{proof}
\end{thm}

\begin{thm}
 Suppose that assumptions of Theorem \ref{thm:Rlin} hold. Then, for every $\varepsilon\in(0,e^{-1})$, there exist $\hat\rho\in(0,1)$ and $Q>0$
 such that Algorithm \ref{alg:MBSM} achieves $f_{\N{}}(x_k)-f_{\N{}}(x_*)<\varepsilon$ in at most $k_\varepsilon$ iterations,
 where
 \begin{equation}
  k_\varepsilon=\left\lceil\frac{\log(f_{N{}}(x_0)-f_{N{}}(x_*)+Q)+1}{|\log(\hat\rho)|}\log(\varepsilon^{-1})+1\right\rceil.
 \end{equation}
\begin{proof}
 We follow the proof of Theorem 6 in \cite{Grapiglia:complexity}.\\
 In Theorem \ref{thm:Rlin} we proved that the following inequality holds
 \begin{equation}\label{rlj}
  \fsub{}{k}-\fsub{}{*}\leq\rho^k\left(\fsub{}{0}-\fsub{}{*}\right)+\sum_{j=0}^{k-1}\rho^j(\zeta_{k-j}+\omega_{k-j})
 \end{equation}
 where $\omega_j=2\nu_j+C\eta_j$.
 Moreover, under our assumptions the last sum converges to 0 $R$-linearly and hence there there exist $\bar\rho\in(0,1)$ and $Q>0$ such that
\begin{equation}
 \sum_{j=0}^{k-1}\rho^j(\zeta_{k-j}+\omega_{k-j})\leq Q\bar\rho^k.
\end{equation}
Replacing this inequality in \eqref{rlj} and denoting with $\hat\rho$ the maximum between $\rho$ and $\bar\rho$, we get
\begin{equation}\label{rlq}
 \fsub{}{k}-\fsub{}{*}\leq\hat\rho^k\left(\fsub{}{0}-\fsub{}{*}+Q\right).
\end{equation}
We denote with $\ke$ the first iteration such that $f_{\N{}}(x_k)-f_{\N{}}(x_*)<\varepsilon$, so that by \eqref{rlq} we have
\begin{equation}
 \varepsilon<f_{\N{}}(x_{\ke-1})-f_{\N{}}(x_*)\leq\hat\rho^{\ke-1}\left( f_{\N{}}(x_0)-f_{\N{}}(x_*)+Q\right)
\end{equation}
and hence,
\begin{equation*}
 \log( f_{\N{}}(x_0)-f_{\N{}}(x_*)+Q)+\log(\varepsilon^{-1})>-(\ke-1)\log(\hat\rho)=(\ke-1)|\log(\hat\rho)|.
\end{equation*}
Rearranging this expression, since $\log(\varepsilon^{-1})>1$ we get
 \begin{equation}
  k_\varepsilon-1<\frac{\log(f_{N{}}(x_0)-f_{N{}}(x_*)+Q)+1}{|\log(\hat\rho)|}\log(\varepsilon^{-1})
 \end{equation}
 and this completes the proof.
\end{proof}

\end{thm}

\section[Numerical Results]{Numerical Results}
In this section we report on the numerical experimentation we carried out over the methods introduced in Section 2.
Our main goals in this section are the following: to asses whether the subsampling procedure provides a reduction in the overall computational cost,
and to provide an indication of which could be the best combination of choices
for the subsampling schedule and the gradient displacement vector. Then, the subsampled methods will also be compared with the
spectral gradient method without subsampling (SGFull).\\

We considered binary classification problems \cite{nocedal:OMLSML}. Given a dataset made of $\hat N$ of pairs $\{(a_j,b_j)\}_{j=1}^{\hat N}$
with $a_j\in\R^n$ and labels $b_j\in\{-1,1\}$,  we use the 95\% of data as training set and the remaining 5\% as validation set. 
Then, letting $\cal{N}$ and $\cal{V}$  be the sets of indices  of the data in the training  and in the validation set, respectively,    
 we define the objective function $f_{\N{}}$ as
 \begin{equation}\label{eq:binaryclass_objfun}
 f_{\N{}}(x)=\frac{1}{N} \sum_{j\in {\cal{T}}} f_j(x), \phantom{spazio} f_j(x)=\log\left(1+\exp(-b_ja_j^tx)\right)+\lambda\|x\|^2
 \end{equation}
 where $N$ is the cardinality of $\cal{N}$ and  the regularization parameter $\lambda$ is set to $N^{-1}.$\\

Note that the main cost in the computation of each function $f_j$ is the evaluation of $\exp(-b_ja_j^tx_j).$
Computing the gradient of $f_j$ for a generic index $j$, we get
\begin{equation}\label{eq:gradexp}
 \nabla f_j(x)=\frac{\exp(-b_ja_j^tx)}{1+\exp(-b_ja_j^tx)}b_ja_j+2\lambda x,
 \end{equation}
thus the evaluation of the gradient comes for free from the evaluation of the corresponding function.
Relying on this remark, at the beginning of every iteration the values $\exp(-b_ja_j^tx_k)$
are computed and then exploited during the execution to evaluate both the sampled function
$\fsub{k}{k}$ and the sampled gradient $\gfsub{k}{k}$. \\

Given an initial size $N_0\geq1$ and a growth factor $\tau>1$, we set $N_k=\lceil\tau^kN_0\rceil$ if this quantity
is smaller than the full sample size $N$, and $N_k=N$ otherwise. We employ \eqref{eq:LiSub} with  $c_1=10^{-4}$ and $
\zeta_k=10^2k^{-1.1}$ and we compute $\alpha_k$ through a backtracking strategy. As a result we have $\alpha_k=\beta^{-\bar \jmath}$
where $\beta=0.5$ and $\bar \jmath$ is the smallest positive integer such that \eqref{eq:LiSub} is satisfied, provided $\bar \jmath\leq 15.$
When such a $\bar \jmath$ does not exist, if the full subsample has not been reached yet, we set $x_{k+1}=x_k$
and we proceed to the next iteration enlarging the subsample. If we are working with the full sample, we declare failure.
We do not explicitly check if Wolfe condition \eqref{eq:WolfeSub} holds at the chosen $\alpha_k$ since the backtracking strategy and the safeguard $\bar \jmath\leq15$ prevent the step size to become too small.\\

We consider three problems of the form \eqref{eq:binaryclass_objfun}, corresponding to the datasets given in
Table \ref{datasets}, along with their dimension $\hat N$ (number of data points  including both training and validation set) and $n$ (dimension of the decision variable).\\

\begin{table}[h]
\begin{center}
\caption{Datasets: number of samples, $\hat N$, problem dimension $n$.}
{\begin{tabular}{c|c|c}
\hline
 \textbf{Dataset}&$\hat N$& $n$\\
 \hline
 CINA0$\ $\cite{set:cina} & 16033 & 132\\
 MNIST$\ $\cite{set:mnist} & 60000 & 784\\
 Mushrooms$\ $\cite{set:mnist} & 5000 & 112\\
 \hline
 \end{tabular}}
\label{datasets}
\end{center}
\end{table}

 \subsection{Comparison of the Subsampled Methods}

 We report the results of the comparison among the performance of all the subsampled methods summarized in Table \ref{tabellametodi}
 on the CINA0 dataset. We underline that the reported statistics are representative of the results that we obtained also with
the other datasets. We add comparison with the method employing exact functions and gradients (SGFull) for completeness.

 We run each method with $\tau=1.1$, $N_0=3$, and we stop the procedure when
 the norm of the gradient is smaller than $10^{-4},$ provided that the full sample is reached.
 We count the number of iterations performed, the number of scalar products computed for the evaluation of each term
 $\exp(-b_ja_j^tx_j) $ with  $j\in {\cal{N}}_k$, the number of functions evaluations and gradients evaluations.
 We recall that since the subsamples are chosen randomly the obtained sequence is not deterministic and therefore for each method
 we perform 100 runs and we report the averages. \\

 In Table \ref{tabellatest1} we report the obtained results. The averages obtained are
 divided by the size $N$ of the training set, so that the evaluation of the full gradient and function counts 1.
 For the number of gradients evaluations, we consider two countings. In column GE\_1 we count only
 the computation of gradients $\nabla f_j(x_k)$ corresponding to functions $f_j$ that have not been previously evaluated at
 the same point $x_k$,
 while in GE\_2 we count every gradient evaluations irrespectively to the evaluation of the functions.
 The motivation behind this choice is the remark that we made at the beginning of this section about the particular
 form of the gradients for the problems we are considering. Since the evaluation of the gradient through equation \eqref{eq:gradexp}
 does not require additional computation with respect to the corresponding evaluation of the function,
 we can count only the gradients evaluations involving new scalar products (that is, the gradients required to
 compute the vector $y_k$); we choose to report also the full count of the gradients to better understand what
 would be the computational cost if we could not rely on this property of the gradients, that is
 in case we consider a sum of functions $f_j$ such that there is not this strict relationship between gradients and functions.\\
 In the table the headers of the columns have the following meaning: IT is the number of iterations performed, SP is the number of
 scalar products computed, FE is the number of functions evaluations and GE\_1 and GE\_2 are the gradients evaluations explained
 above.\\
 \begin{center}
 \begin{table}[H]
 \begin{center}
 \caption{ Statistics of the runs}
 { \begin{tabular}{c|c|c|c|c|c}
   \hline
   METHOD & IT & SP & FE & GE\_1 & GE\_2 \\
   \hline
   SG\_N\_1 & 101 & 67.6 & 66.5 & 1.1 & 31.3\\
   SG\_N\_2 & 106 & 80.1 & 80.1 & 0 & 35.9\\
   SG\_I\_1 & 109 & 91.6 & 85.7 & 5.7& 44.7\\
   SG\_I\_3 & 105 & 93.6 & 93.6 & 0 & 34.7\\
   SGFull & 52& 115 & 115 & 0 & 52\\
   \hline
  \end{tabular}}
\label{tabellatest1}
\end{center}
\end{table}
\end{center}
We can notice from Table \ref{tabellatest1} that the nested subsampling methods are in general more efficient than their non-nested counterparts
and that the definition of $y_k$ that uses the same subset $\N{k}$ to approximate the gradient both at $x_k$ and at $x_{k-1}$ seems to
perform better than the alternatives, both in the nested and in the non-nested case.
In fact, our numerical experience suggests that computing the gradient displacement vector using all the current information, i.e.
using the gradient sampled in the same set $\N{k}$ where the function and the search direction are sampled, allows to obtain a less
expensive method in terms of weighted number of functions evaluations. Then, even if choices \eqref{eq:Y2} and \eqref{eq:Y3}
do not require extra scalar products for computing $y_{k-1},$ the gain obtained is not enough to compensate the larger number of
functions evaluations required. Finally, notice that the results indicate that using the subsampled approach is overall more efficient than SGFull.

 \subsection{Influence of the Parameter $\tau$}
In this subsection we focus on the role of the growth factor $\tau$. The influence of
this parameter over the efficiency of the methods may be high:
a bigger value of $\tau$ means a more accurate approximation of the objective function and its gradient from the beginning of
the algorithm, but it also causes a higher per-iteration cost during the initial phase. We consider the methods SG\_N\_1 and SG\_I\_1 which from the previous section appear to be the best among the nested and the non-nested methods,
respectively.

In Figure \ref{tau_test} we
report the number of scalar products required by the two methods for  different values of $\tau \in \{1.1,1.2,...,1.9, 2,2.25,2.5,...,4.5, 5\}$.
For each $\tau$ we consider
100 runs of each method on the CINA0 dataset and we take the average number of scalar products computed excluding the worst 20  and the best 20  results,
then we divide by $N$.
We compare the overall cost in terms of scalar products with that of  SGFull.\\
In Figure \ref{tau_test_var} we plot the difference between the smallest and the biggest number of scalar products required among
the runs of each of the two methods: this can be taken as a first indication of the variance in the efficiency of the methods
corresponding to the chosen $\tau$.\\

 \begin{minipage}{0.45\textwidth}
\begin{figure}[H]
\centering
\includegraphics[width=0.95\textwidth,height=100pt]{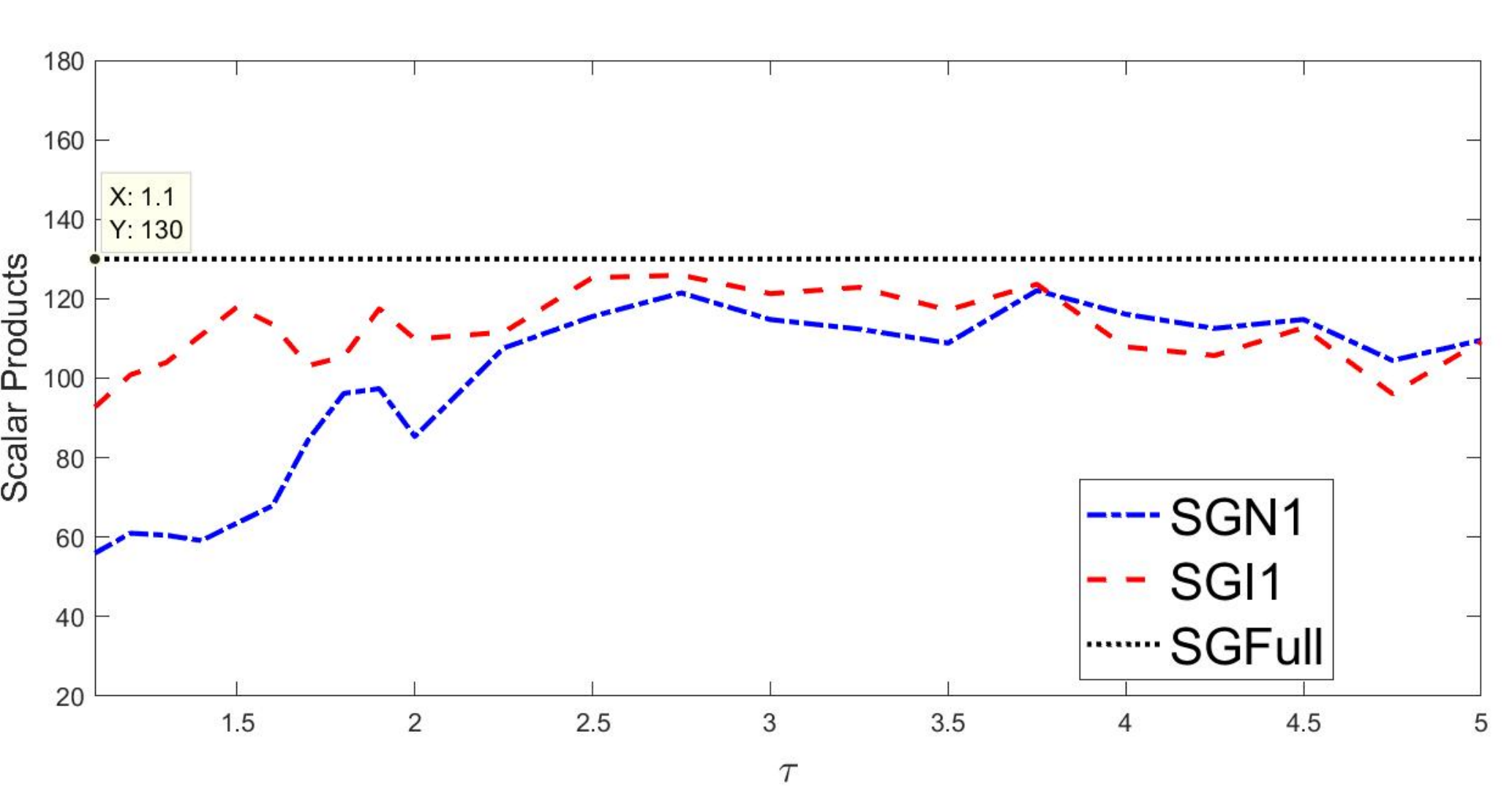}
\caption{}\label{tau_test}
\end{figure}
\end{minipage}
\begin{minipage}{0.45\textwidth}
\begin{figure}[H]
\centering
\includegraphics[width=0.90\textwidth, height=100pt]{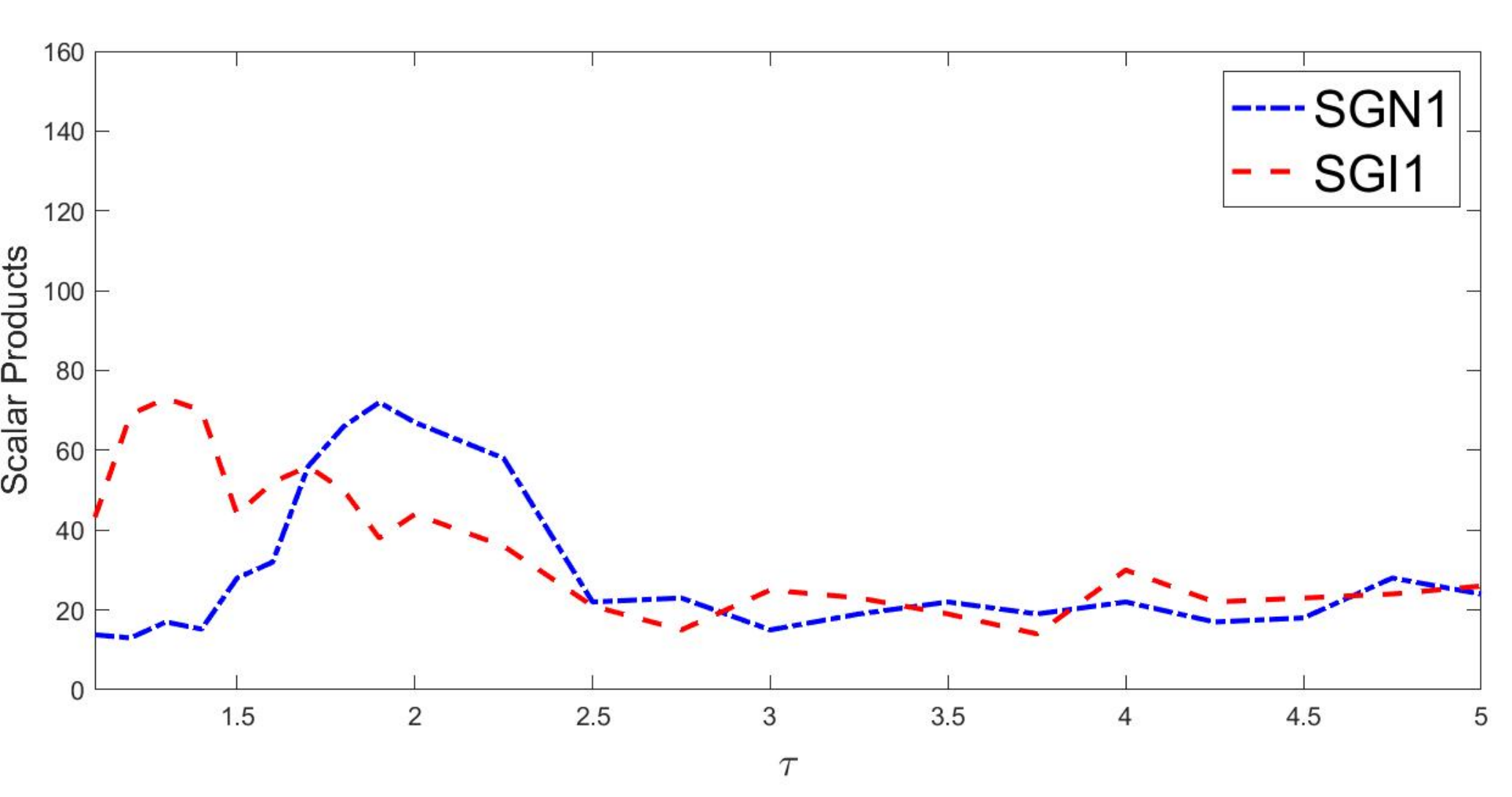}
\caption{}\label{tau_test_var}
\end{figure}
 \end{minipage}\\
 $\ $\\
As we can see from Figure \ref{tau_test}, among the considered values of $\tau$ the initial
choice $\tau=1.1$ seems to be the optimal one for both the methods and, for small values of the growth factor SG\_N\_1 seems to be more
efficient than SG\_I\_1.
As $\tau$ grows the computational cost of the nested method tends to increase (for $\tau>3.5$ it becomes higher than the cost
of SG\_I\_1). For all the tested values of $\tau$ the subsampled methods outperform the SG method with full function and gradient information.
The gap is more evident for small values of $\tau.$ For larger values the benefits deriving from a more accurate representation of $f_{\N{}}$ do not seem to
be enough to compensate for the bigger average number of per-iteration
scalar products required.
Figure \ref{tau_test_var} shows that the cost of the non-nested method seems to be overall less influenced by the choice of the
growth factor.   We considered also values of $\tau$ smaller than 1.1 and we observed that they give rise to procedures overall more expensive than that corresponding to $\tau=1.1$ as the growth of the data set is too low   and starting from $N_0=3$ the number of iterations needed to reach  a reasonable value of the training loss is too high.

 \subsection{Training Error}\label{trainingerror}
 We consider the three datasets reported in Table \ref{datasets} and three methods:
 SGFull, SG\_N\_1 and SG\_I\_1 both with with $\tau=1.1$.
 For every method and every dataset we study how the \emph{training error} $f_{\N{}}(x_k)-f_{\N{}}^*$ changes as the number of
 iterations and  scalar products performed grows. We stress that the number of performed scalar products can be considered a
 measure of the overall computational cost due to the form of functions $f_j$ and their derivatives.
 The optimal value $f_{\N{}}^*$ has been computed running SGFull with termination condition $\|\nabla f_{\N{}}(x_k)\|\leq10^{-7}$.\\
  \vspace{-15px}
  $\ $\\
\begin{figure}
\includegraphics[width= 0.45\textwidth, height=105pt]{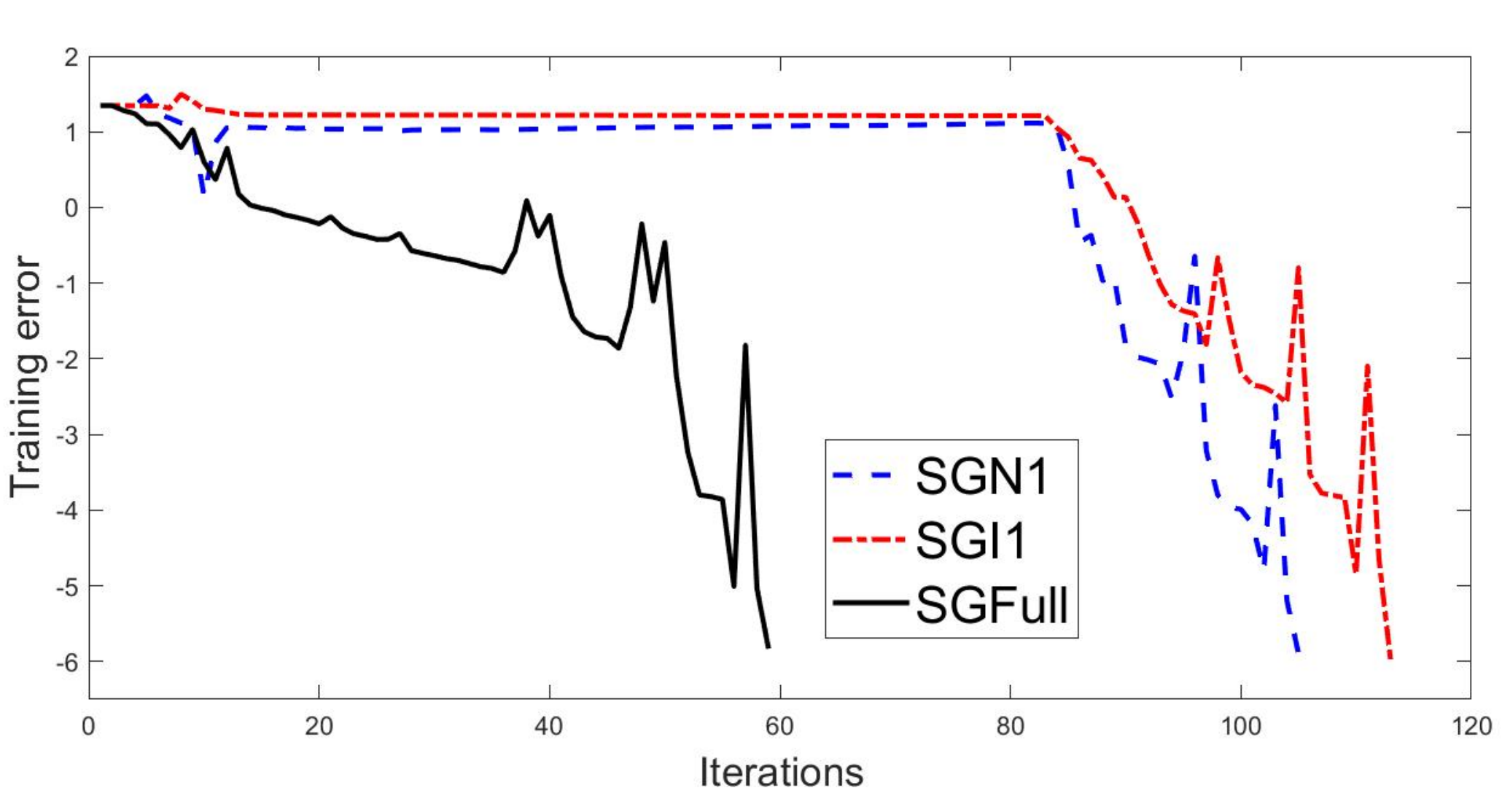}
\includegraphics[width=0.45 \textwidth, height=100pt]{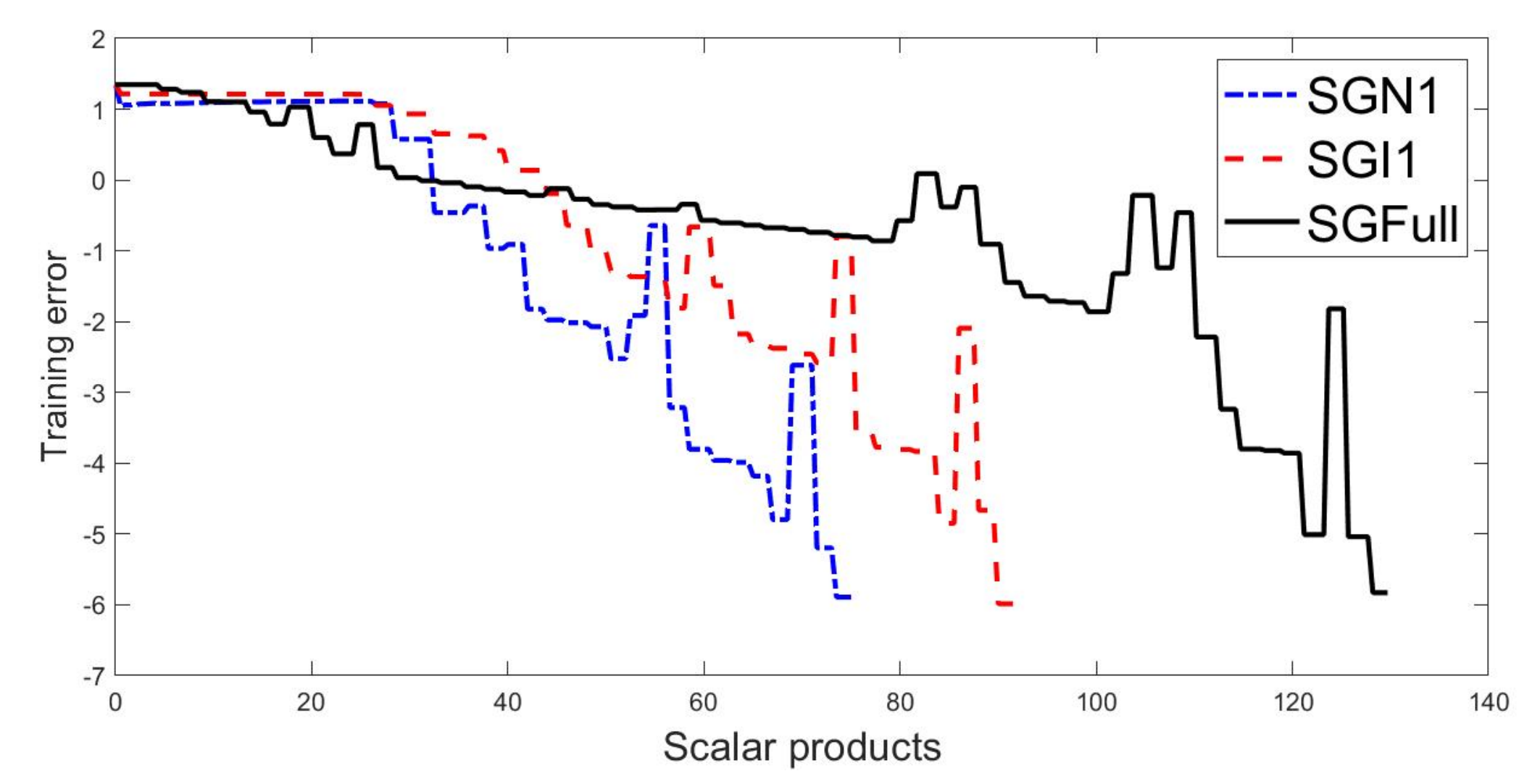}
\caption{CINA0 dataset, training error versus iterations (left) and versus scalar products (right).}\label{psVSerr_cina}
\end{figure}
\begin{figure}
\includegraphics[width=0.45\textwidth,height=100pt]{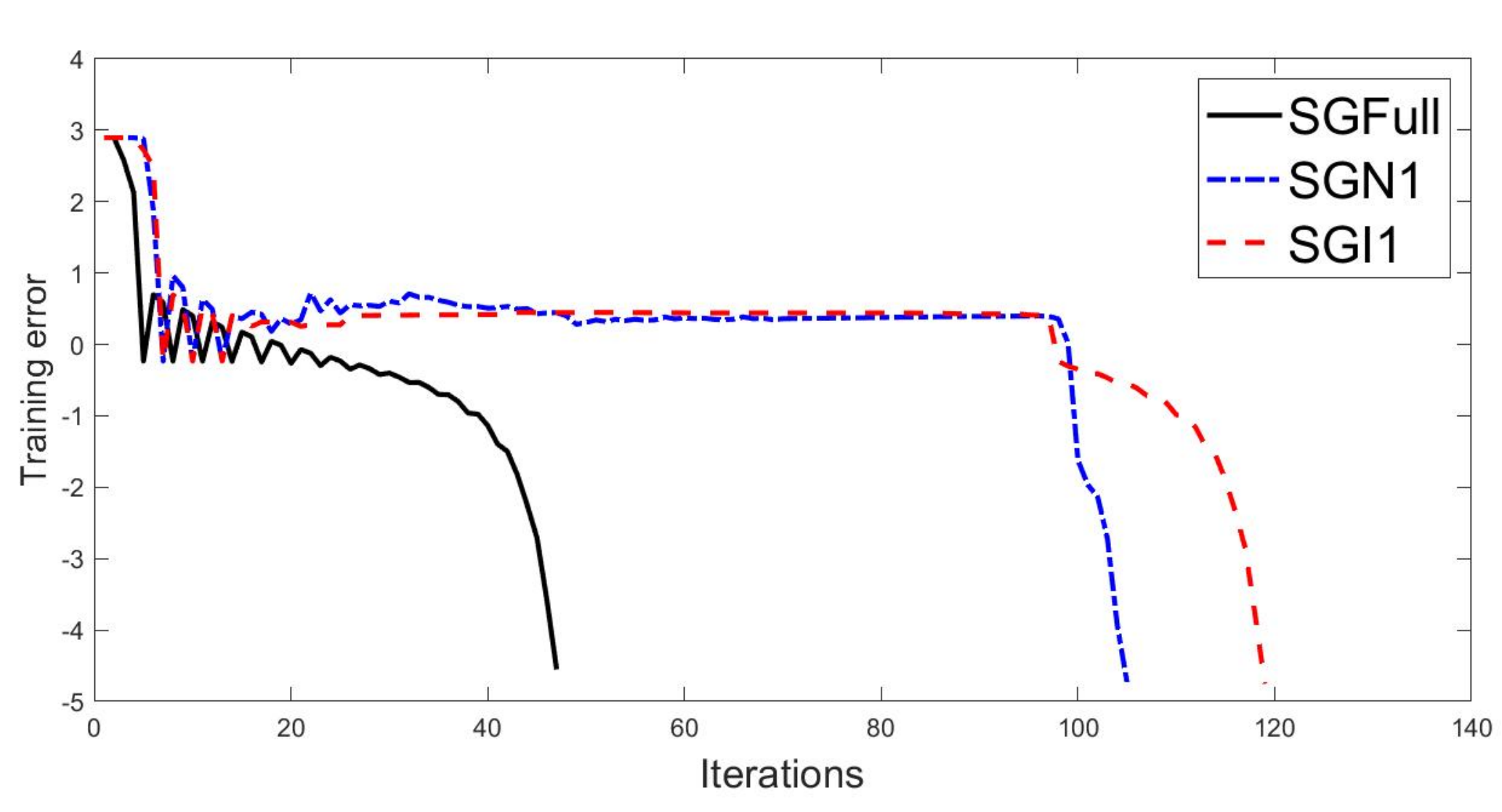}
\includegraphics[width=0.45\textwidth, height=100pt]{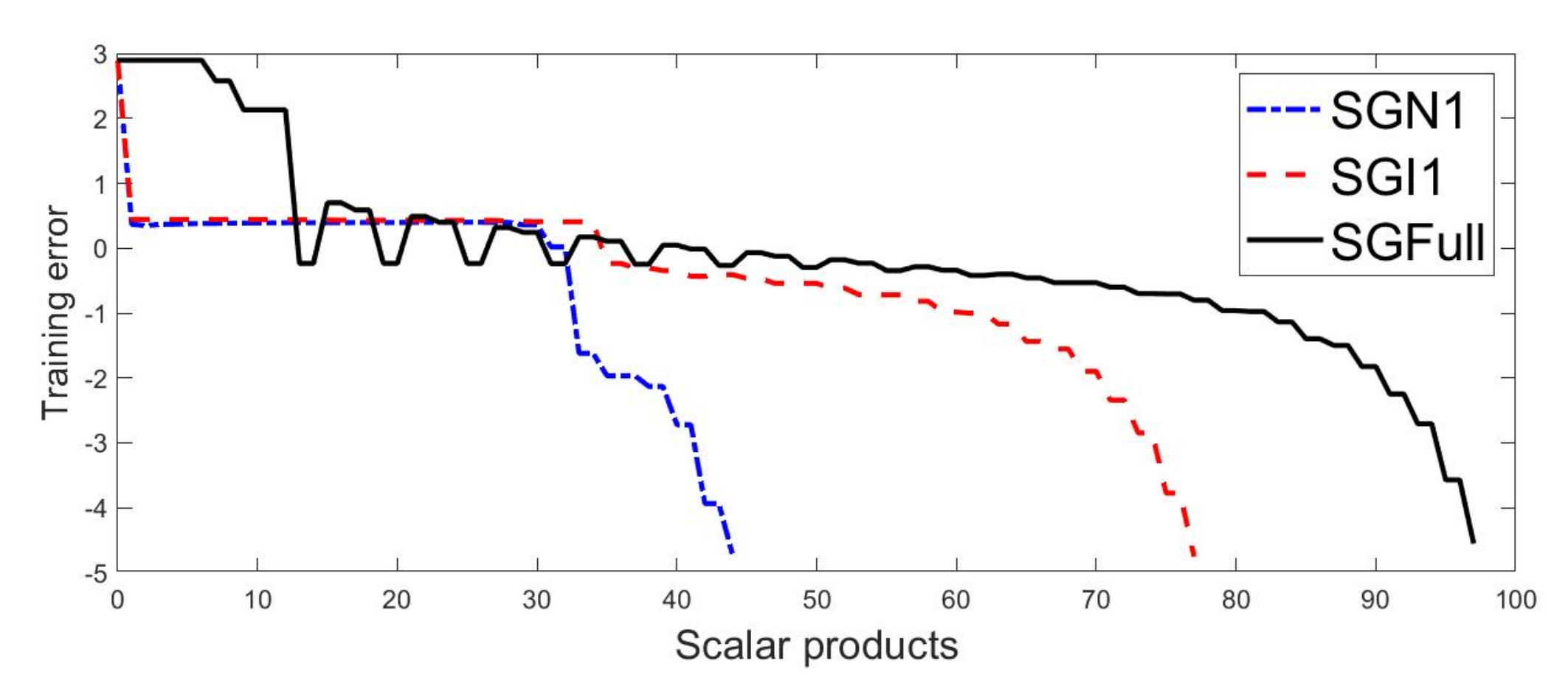}
\caption{MNIST dataset, training error versus iterations (left) and versus scalar products (right).}\label{psVSerr_mnist}
\end{figure}
\begin{figure}
\includegraphics[width=0.45\textwidth,height=100pt]{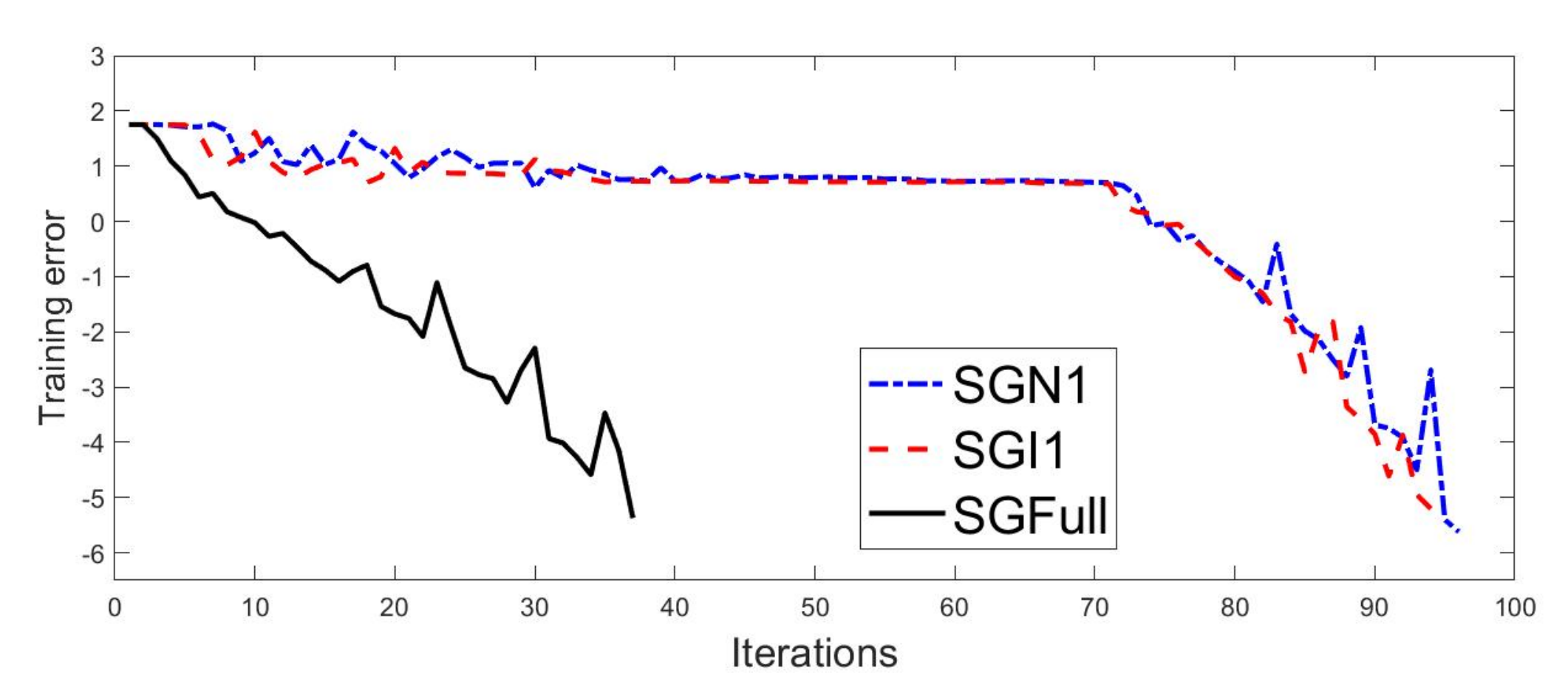}
\includegraphics[width=0.45\textwidth, height=100pt]{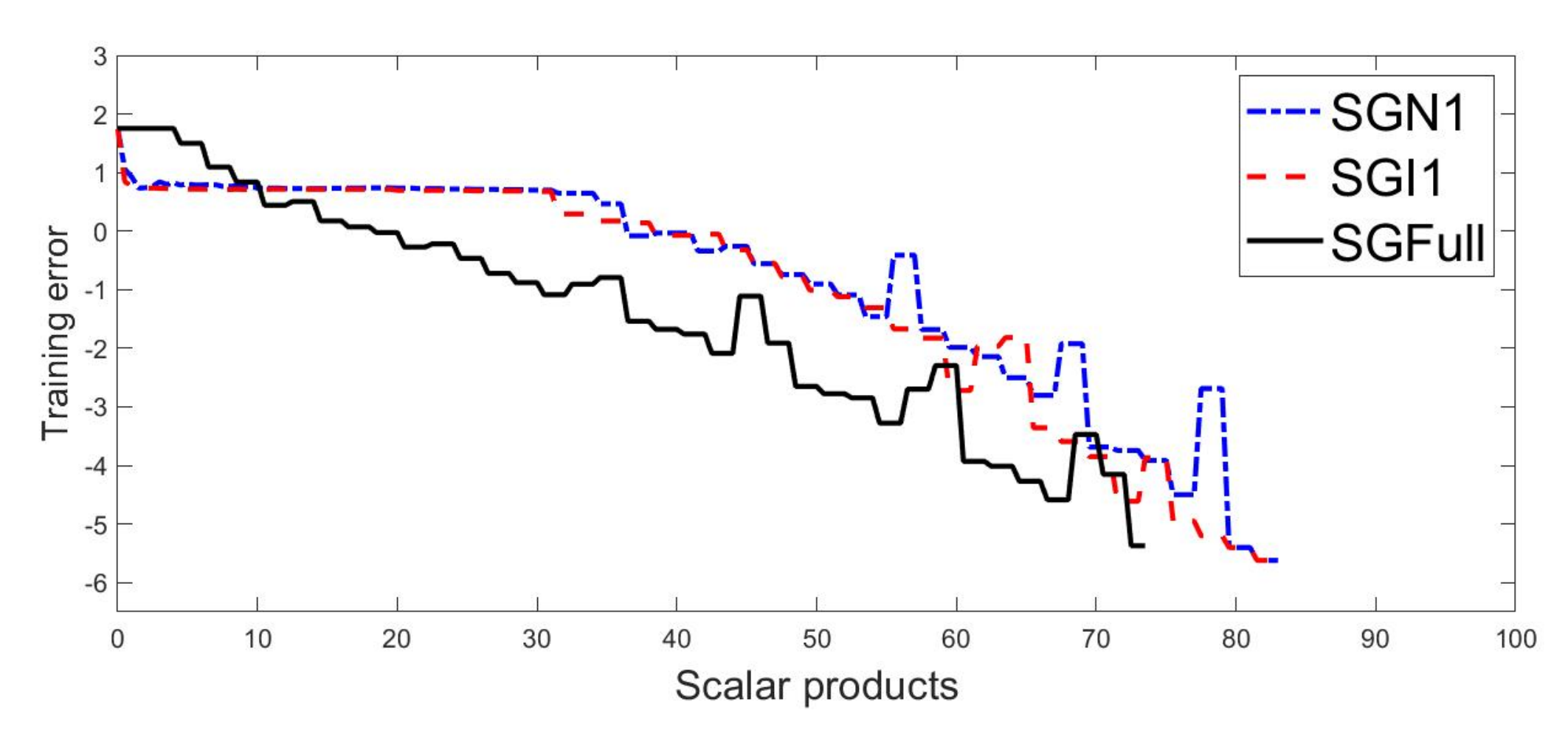}
\caption{Mushrooms dataset, training error versus iterations (left) and versus scalar products (right).}\label{psVSerr_mush}
\end{figure}
In Figures  \ref{psVSerr_cina}-\ref{psVSerr_mush} we plot the training loss versus iterations on the left and versus scalar products on the right. The reported figures refer to a specific run out of $100$ runs,  representative of the average behavior of the methods. Plots are in logarithmic scale on the $y$-axes.
As expected, for all the three datasets, the spectral gradient method without subsampling requires a much smaller number
 of iterations with respect to SG\_N\_1 and SG\_I\_1.
  Let us now consider the number of scalar products.  Concerning CINA0 and MNIST dataset, both the subsampled methods are less expensive than SGFull, that is, they produce the same training error with a smaller number of scalar products. Moreover, the nested method
 SG\_N\_1 appears to be better than the non-nested one.
 Regarding the Mushrooms dataset, we have that SG\_N\_1 and SG\_I\_1 appear to behave very similarly and they are both less efficient
 than SGFull, therefore on this particular dataset the subsampling scheme does not seem to be convenient with respect to the method
 without subsampling.
 This is due to the fact that the size of the training set is small and therefore the gain obtained reducing the sample size is not
 enough to produce an overall saving.

 We also notice that we  do not expect an advantage in using
  the subsampling in all the situations where problems are easy enough   (because the starting point is close to
 the solution, or due to the particular form of the functions $f_j$) to be solved by
SGFull within  a small number of iterations.

\subsection{Validation Loss}

In the previous subsections we tested the behavior of the Algorithm 2.2 using the deterministic stopping criterion - the norm of the full gradient. Then, we always reach the full sample
and ask for an $\epsilon$-accurate first-order method. However,  in these applications  a reasonable value of the  validation error is needed   rather than an   accurate approximation of the minimizer of  \eqref{eq:binaryclass_objfun}.

In order to do that, we show numerical  results obtained using a  stopping criterion related to the validation loss  rather than to the training loss or to the full gradient. Notice that the testing data set is usually much smaller than the training set, in our runs it is constituted by the $5$\% of samples of the whole dataset.
Letting 
$$
 f_{\V{}}(x)=\frac{1}{|\V|} \sum_{j\in {\cal{V}}} f_j(x), 
 $$
be the validation loss, 
we stop whenever the following condition is met:
\begin{equation}\label{val_stop}
f_{\V}(x_{k})>1.1 f_{\V}(x_{k-1}) \mbox{ or  } | f_{\V}(x_{k-1})-f_{\V}(x_k))|<10^{-3} | f_{\V}(x_{k})) |
\end{equation}
provided that $N_k\ge p N$, with $p \in (0,1)$. Thus, we stop whenever we monitor a 10\% increase or a stallation  of the validation loss,  provided that the sample size is at least a fixed percentage of the full sample size.

\begin{figure}
\centering
\includegraphics[width=0.45\textwidth,height=100pt]{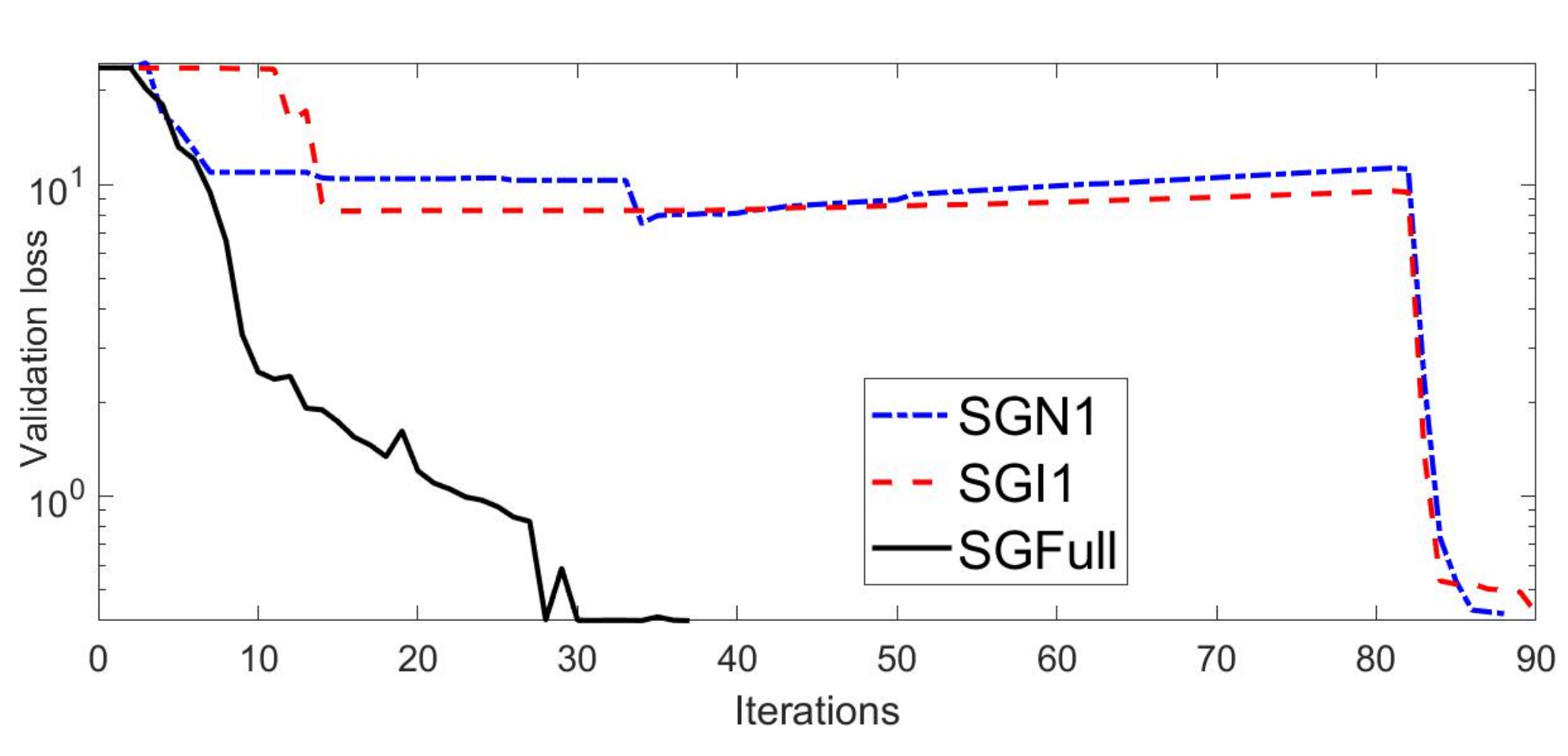}
\includegraphics[width=0.45\textwidth,height=100pt]{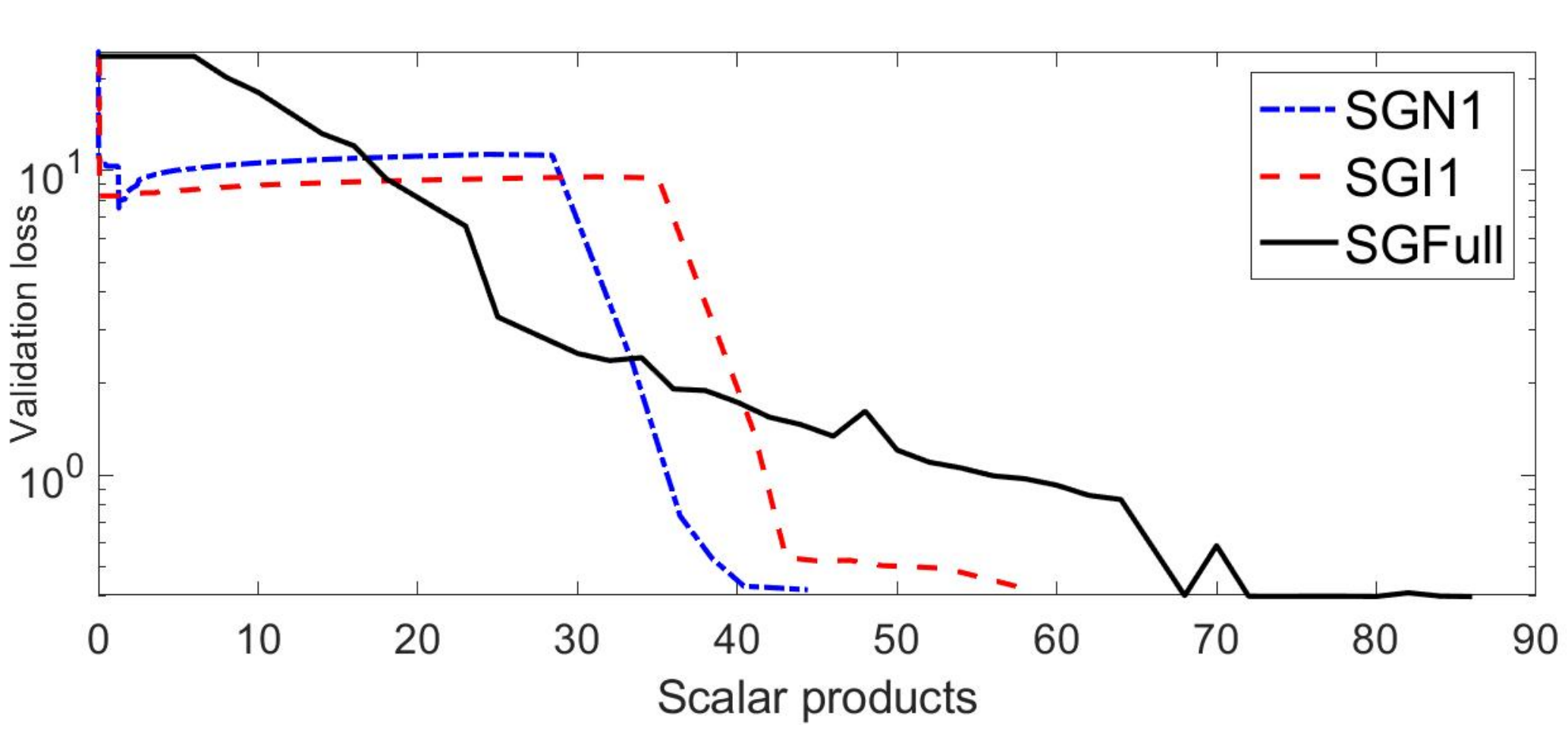}
\caption{CINA0 dataset, validation loss versus iterations (left) and versus scalar products (right),  stopping criterion \eqref{val_stop} with $p=1$}\label{Cinaval}
\end{figure}

\begin{figure}
\centering
\includegraphics[width=0.45\textwidth, height=100pt]{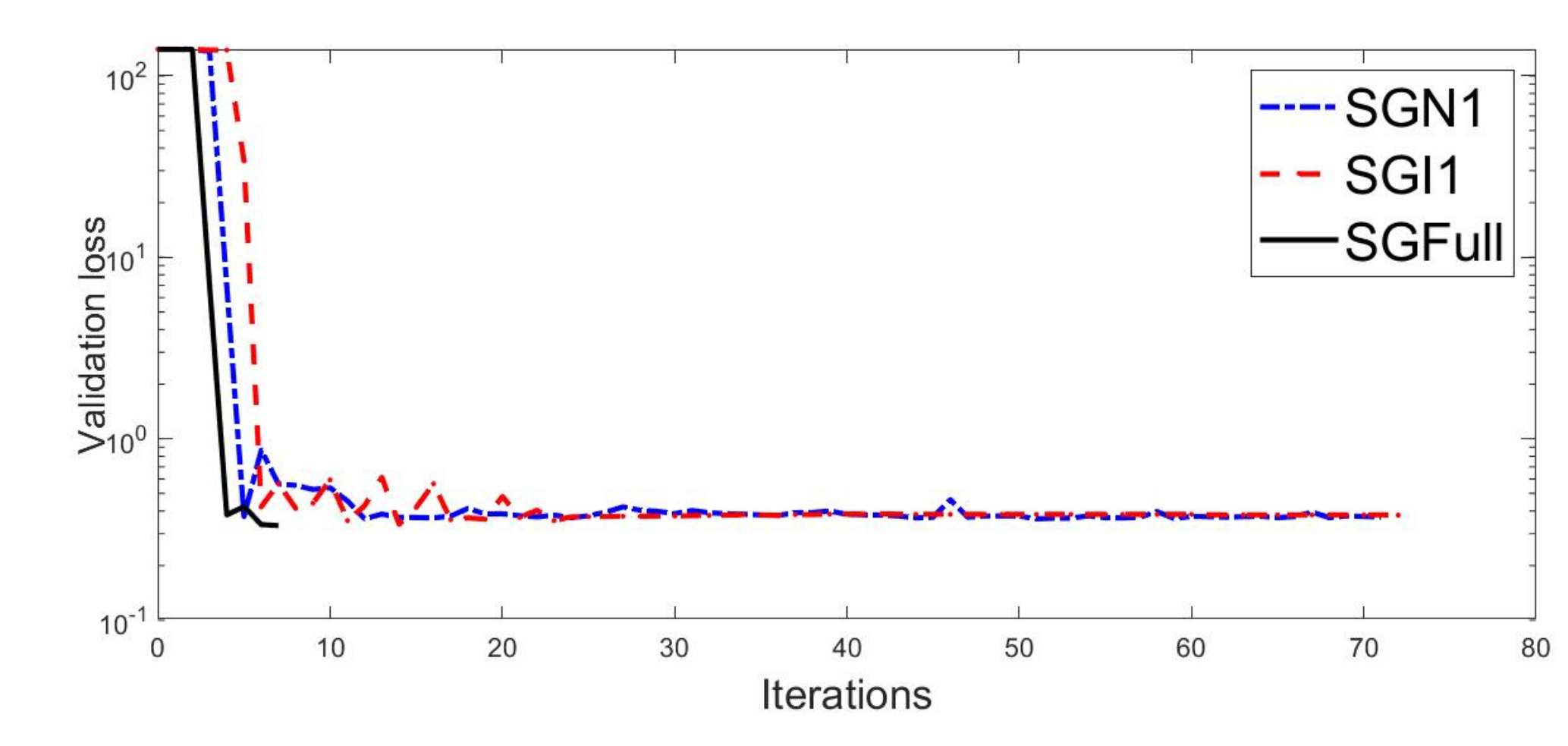}
\includegraphics[width=0.45\textwidth, height=100pt]{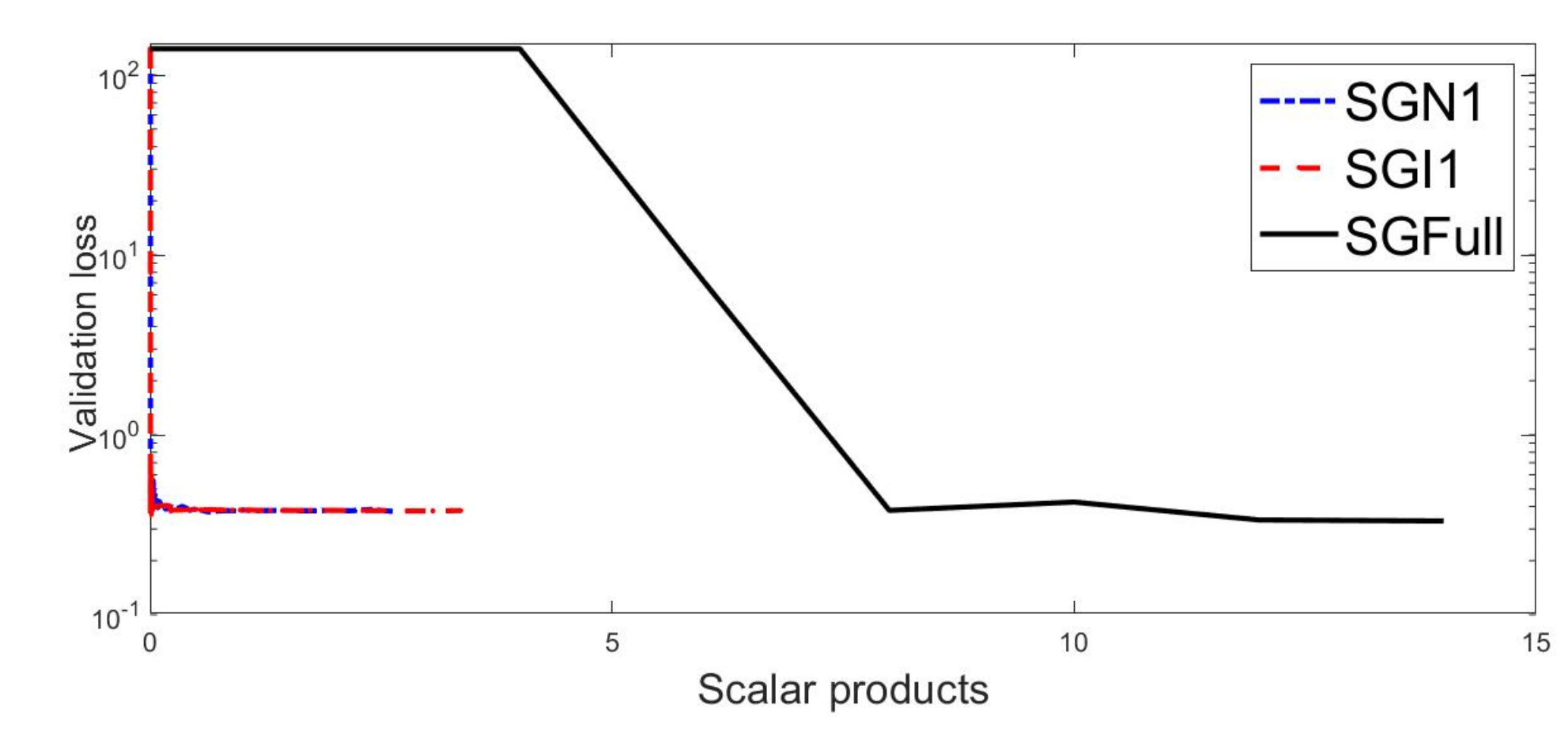}
\caption{MNIST dataset, validation loss versus iterations (left) and versus scalar products (right),   stopping criterion \eqref{val_stop} with $p=0.1$}\label{MNISTval}
\end{figure}

In Figures \ref{Cinaval}-\ref{MNISTval} we plot the validation loss versus iterations and scalar products for both SG\_N\_1 and SG\_I\_1 using logarithmic scale on the $y$-axes.
For the sake of comparison we also plot the validation
loss  of SGFull stopped  when
 the norm of the gradient is smaller than $10^{-4}$.   This way we give a clear indication of the value of the validation loss that can be obtained by a method computing a $10^{-4}$-accurate first-order critical point.  The reported figures refer to a specific run out of $100$ runs,  representative of the general behavior of the methods.
 As a first comment we observe that in case of CINA0 dataset, we needed to reach the full sample in order to provide a validation error of the same order of that provided by SGFull. Therefore, for CINA0 dataset, we use  $p=1$ in \eqref{val_stop}. In fact, we can observe that the validation decreases in the very early stage of the iterative process, then it remains almost constant, and then rapidly drops to $0.4$ as soon as the full sample is reached. However, notice that only a very few  iterations are performed using the full sample as the average number of iterations is $86$ for SG\_N\_1 and $90$ for SG\_I\_1. Notice that the full sample is reached at iteration $82$. On the other hand, considering Mushrooms and MNIST datasets  a smaller number of samples are enough  to provide an acceptable validation error. Thus, we set $p=0.1$ in  \eqref{val_stop}. In case of Mushrooms  the average sample size at termination is $1164$ and $1234$ for SG\_N\_1 and SG\_I\_1, respectively. Both methods used around 25\% of samples and reached  a validation loss value of the order of $0.6$ while SGFull reaches a value of about $0.3$. Analogous behaviour can be  observed in the solution of MNIST problem where  the  average sample size at termination is $7179$ and $6520$ for SG\_N\_1 and SG\_I\_1, respectively. Then, the 10\% of samples are enough for this dataset.

\begin{figure}
\centering
\includegraphics[width=0.45\textwidth,height=100pt]{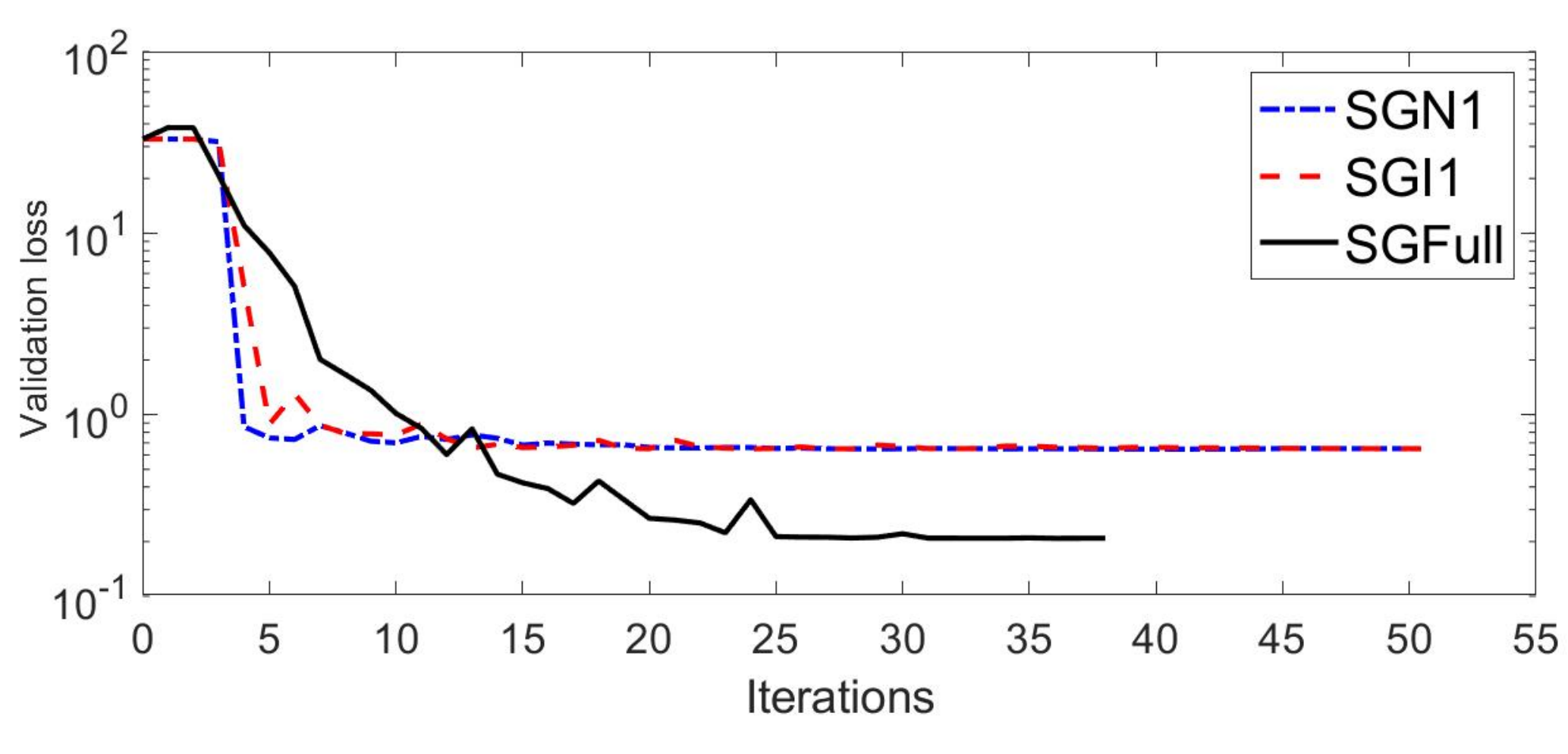}
\includegraphics[width=0.45\textwidth, height=100pt]{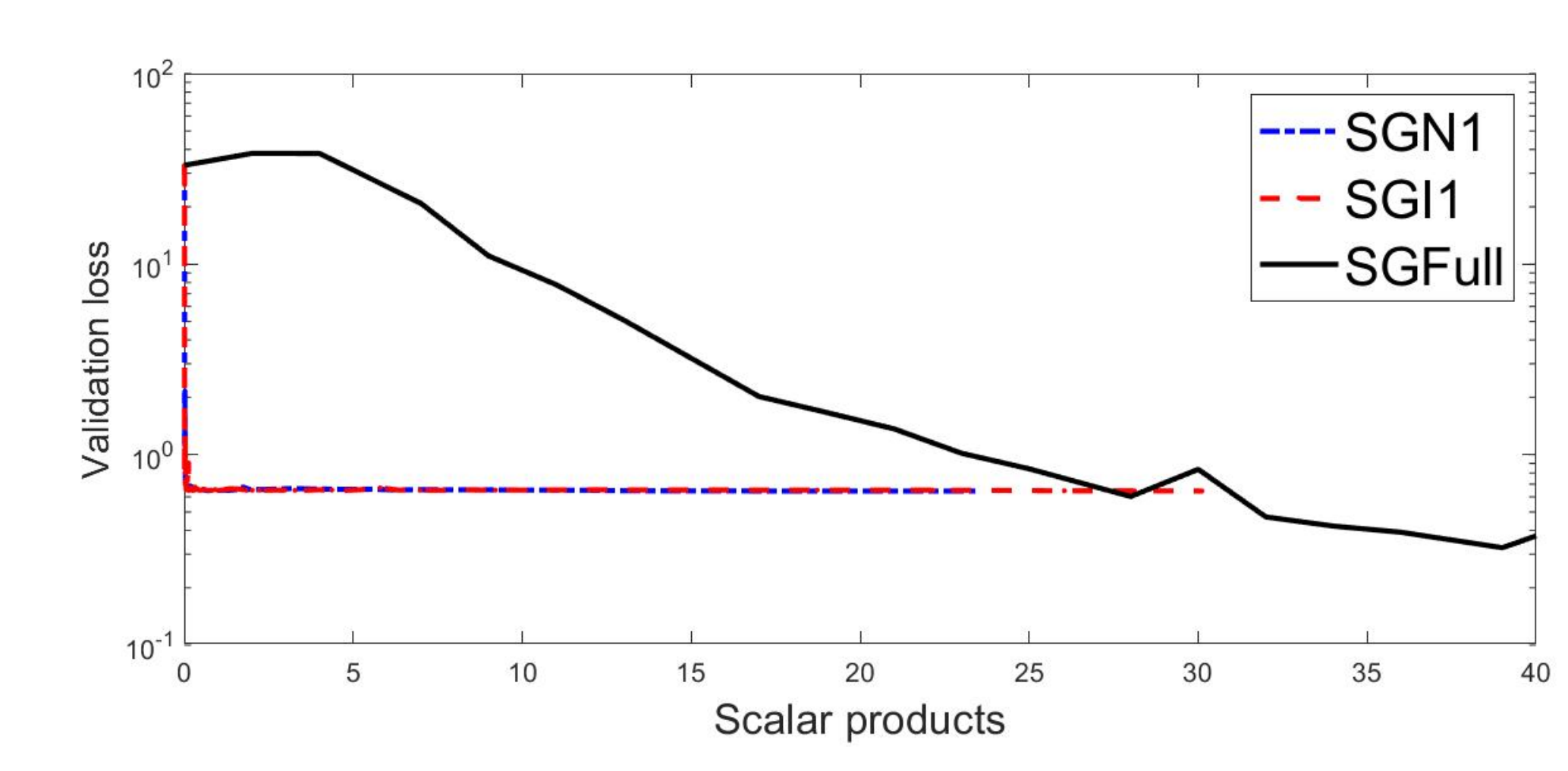}
\caption{Mushrooms dataset, validation loss versus iterations (left) and versus scalar products (right),  stopping criterion \eqref{val_stop} with $p=0.1$}\label{Mushval}
\end{figure}

\section{Conclusions}

We analyzed subsampled spectral gradient methods for solving large sum optimization problems with continuously-differentiable objective function. Our main aim was to provide initial understanding of the behavior of such methods with different kinds of spectral coefficients, or more precisely, with different kinds of gradient difference calculations  used for obtaining spectral coefficients within the growing sample size framework. 

 Although the safeguards for the spectral coefficients provide a descent direction regarding the current approximate function, globalization strategy relies on nonmonotone  line search. The motivation for this is twofold - to the best of our knowledge, both spectral gradient methods and subsampling methods favor  nonmonotone line search techniques.  Moreover, if the objective function is assumed to be  strongly convex, R-linear convergence is achieved. We provided complexity results for both convex and non-convex case.  Further modifications of this approach could be devised in order to solve optimization problems arising in the training of neural network employing non-smooth activation functions.

According to the tests performed on three binary classification problems, our main conclusions are as the following: 1) nested (cumulative) samples perform better than non-nested ones considering costs of the algorithm in terms of scalar products; 2) employing the same subset ${\cal N}_k$ for computing the displacement  vector $y_k$ proved to be beneficial despite the needed 
additional cost;  3) subsampling seems to reduce the overall computational cost  as the tested subsampling variants outperform the full spectral gradient methods 
except in Mushrooms dataset where the size of the training set is small.

Finally, investigating on the validation loss, we conclude that a mini-batch version of the proposed spectral gradient methods has a good potential. However, stopping criterion in this setting is a problem itself and it needs further analysis. A related question is how the proposed methods compete in the unbounded sample case. This will be one of the topics of our further research. 

\section*{Acknowledgements}
\begin{small} Indam-GNCS partially supported the first author  under Progetti di Ricerca 2018. The second author was supported by the
Serbian Ministry of Education, Science and Technological Development, grant no. 174030. 
 The third author was supported by the European Union's Horizon 2020 programme under the Marie Sk\l odowska-Curie Grant Agreement No 812912.
The authors are grateful to the anonymous referee whose suggestions led to significant
improvement of the manuscript. \end{small}
\bibliographystyle{plain}
\bibliography{bibliography}

\end{document}